\def\rpp{\RR P^2} \def\lat{\rm}
\input amstex
\documentstyle{amsppt}
\input tmm.def
\input virodef.tex
\hsize=348pt
\vsize=572pt
\def\ra{\RR A}
\def\ca{\CC A}
\def\R{\Bbb R}
\NoBlackBoxes
\topmatter
\at  Complex orientations of real algebraic surfaces \endat
\au Oleg Viro \endau
\af
Mathematics Department, 
Uppsala University,
Uppsala, 752 06, 
Sweden;\newline
POMI, Fontanka 27,
Sankt Petersburg  191011,
Russia.
\endaf
\toc
\ah 0. Introduction \endah
\ah 1. Generalization of type of a curve  \endah
\bh 1.1 Homological reformulation of definition of type\endbh
\bh 1.2 Relation to the form of the complex conjugation
        involution\endbh
\bh 1.3 No more relations among homology classes of real
        components\endbh
\ah 2. Complex orientation of a real surface bounding in
        complexification \endah
\bh 2.1 Digression: Arnold's principle\endbh
\bh 2.2 Complexification of the notion of boundary does
        work, suggesting definition of complex orientation of
        surface\endbh
\bh 2.3 Kharlamov's congruence\endbh
\bh 2.4 Complex orientations and classes lifted from
        the orbit space of complex conjugation\endbh
\ah 3. Relative complex orientations of a real surface  \endah
\bh 3.1 Types of real algebraic surfaces revised\endbh
\bh 3.2 Complex orientations of a real surface modulo
        a curve\endbh
\bh 3.3 Semi-orientations\endbh
\bh 3.4 Internal definition of complex orientations
        (without 4-fold covering)\endbh
\bh 3.5 Orientations modulo changing curve\endbh
\bh 3.6 Orientations of the associated covering surface\endbh
\ah 4. Survey of subsequent results \endah
\bh 4.1 Complex orientations of high-dimensional
        varieties\endbh
\bh 4.2 Spin-structure of a real algebraic surface
bounding in complexification \endbh
\bh 4.3 $Pin^-$-structure of a real algebraic surface of
type $I_{rel}$\endbh
\endtoc
\endtopmatter \document

\ah Introduction\endah

In this paper we will deal with real algebraic {\it surfaces\/}. However let
me start with some well known facts on real algebraic {\it curves,\/}
to motivate the forthcoming exposition.

The set $\ra$ of real points of a non-singular real algebraic curve $A$
lies in the set $\ca$ of its complex points.
Topology of the placement of $\ra$ in $\ca$ can be
described in terms of at most two very concise characteristics.

First,  either $\ra$
divides $\ca$ into two halves, or $\ra$ does not divide $\ca$. In the first
case the curve $A$ is said to be of {\it type\/} I, and called {\it
dividing curve,\/} in the second it is said to be of {\it type\/} II, and
called {\it non-dividing curve.\/}

Second, if the curve is of type I, the canonical orientation of $\ca$
(determined by the complex structure of $\ca$)
determines orientations of both of its halves and they in turn induce
orientations on $\ra$, as on their common boundary. These orientations
are opposite to each other. They are called {\it complex
orientations\/} of $A$.

The type of curve enhanced, in the case of type I, by the complex
orientations describe the inclusion $\ra\to\ca$ up to homeomorphism of
$\ca$.

It was F.~Klein who introduced (more than hundred years ago) the
two types of real algebraic curves, see \cite{K}. In seventies V.~A.~Rokhlin
\cite{R2}, \cite{R3} introduced the complex orientations.
For the case of non-singular real plane projective curves Rokhlin found also
relationships between placement of $\ra$ in $\ca$ and placement of the same
set $\ra$ in the real projective plane $\rpp$. These relationships proved
to be very important in the subsequent development of topology of plane
real algebraic curves.

A traditional viewpoint on problems of topology of real algebraic curves
was that the main problem is to classify up to homeomorphisms of $\rpp$
placements $\ra\subset\rpp$ for non-singular real algebraic curves $A$
of a given degree. At first, the relationships found by Rokhlin seemed
to be useless from this traditional viewpoint. They relate placement of
$\ra$ in $\ca$ to placement of $\ra$ in $\rpp$, but at that time (middle
seventies) they gave nothing new on placement of $\ra$ in $\rpp$ itself.
However later, when additional restrictions were found (the most
remarkable one was found by T.~Fiedler \cite{F}), they gave new
restrictions on placement in $\rpp$ of the real point set of a
non-singular real algebraic curve of a given degree.

Gradually, this development made specialists change the traditional
viewpoint.  The problem of topological classification of placement
$\ra\subset\rpp$ for non-singular real algebraic curves $A$ of a given
degree, as the main problem of the topology of real plane algebraic
curves, has been replaced by a finer problem of topological
classification of placements $\ra\subset\rpp$ taking into consideration
not only degree of $A$, but also placement of $\ra$ in $\ca$.

 Non-singular projective hypersurfaces can be considered as the most
straightforward generalization of non-singular plane projective curves,
and usually results on topology of non-singular real plane projective
algebraic curves have more or less straightforward generalizations to
the case of non-singular real projective hypersurfaces. Surfaces of
the three-dimensional projective space are hypersurfaces next to plane
curves.

Thus, it was natural to expect that the types of curves and the complex
orientations of dividing curves are generalized to high-dimensional
non-singular algebraic varieties, and, first of all, to non-singular
real algebraic surfaces.

This problem was suggested to me by Rokhlin in late seventies, and I
managed to find an answer. However, I could not generalize the most
impressive applications of complex orientations of curves. Therefore I
delayed a detailed publications, restricting myself to short
announcements \cite{V2}, \cite{V3}, \cite{V4}.

Now the subject attracts new people, see \cite{Ka}, \cite{D}. In this
volume several papers are devoted to or motivated by it. The original
definitions and constructions almost are not mentioned there. Almost
the same constructions look much more sophisticated. This is a natural
process, but I do not like that the original faces of the subject have
not appeared in literature.  Therefore I decided to present the
original approach with more or less complete motivations.

\ah 1. Generalization of type of a curve\endah

\bh 1.1 Homological reformulation of the definition of type\endbh
How to reformulate the condition that
$$
\foldedtext\foldedwidth{4.3in}
{the real point set $\ra$ of a non-singular curve $A$ divides the complex
point set $\ca$ of $A$}
$$
in such a way that this condition would have a sense in high-dimensional case?
The codimension of $\ra$ in $\ca$ is equal to the dimension of $A$. Thus if
the dimension of $A$ is greater than 1, and $\ra$ can not divide
$\ca$.

Well, in 1-dimensional case $\ra$ divides $\ca$, if and only if $\ra$ is
zero-homologous in $\ca$, i.~e. it realizes the trivial element of the
group $H_1(\ca;\,\ZZ_2)$. This suggests the following definition:

A non-singular $n$-dimensional real algebraic variety $A$ is said to
{\it bound in complexification\/}, if the set $\ra$ of its real points is
zero-homologous in the set $\ca$ of complex points of $A$,
i.~e. realizes the trivial element of the group $H_n(\ca;\,\ZZ_2)$.

\bh 1.2 Relation to the form of the complex conjugation involution\endbh
Remind that if $\Gt$ is an involution acting in an orientable manifold
$X$ of even dimension $2n$, then by the {\it form of\/}
$\Gt$ one calls the bilinear form
$$
H_n(X)\times H_n(X)\to\ZZ: (\Ga,\Gb)\mapsto \Ga\circ\Gt_*\Gb
$$
where $\circ$ denotes intersection number. This form is symmetric, if
either $n$ is even and $\Gt$ preserves orientation of $X$, or $n$ is
odd and $\Gt$ reverses orientation of $X$. Otherwise it is skew-symmetric.
The inclusion
$H_n(X)\otimes\ZZ_2\to H_n(X;\,\ZZ_2)$ induces isometrical imbedding of
its reduction modulo 2 into the similar {\it $\ZZ_2$-form of $\Gt$\/}
defined by
$$
H_n(X;\,\ZZ_2)\times H_n(X;\,\ZZ_2)\to\ZZ_2:(\Ga,\Gb)\mapsto \Ga\circ\Gt_*\Gb
$$

\tm{1.2.A Lemma {\rm(cf. Arnold \cite{A})}}  Under the
conditions above, if dimension of each component of the fixed point set
$F$ of $\Gt$ is at most $n$, then the union of $n$-dimensional
components of $F$ realizes the characteristic class of the $\ZZ_2$-form
of $\Gt$.  \endtm

(Remind that the characteristic class of a symmetric bilinear form
$b:V\times V\to\ZZ_2$ on $\ZZ_2$-vector space $V$ is a vector $\chi\in
V$ such that $b(\xi,\xi)=b(\chi,\xi)$ for any $\xi\in V$. Any
non-degenerate symmetric bilinear form has a unique characteristic
class. This class is zero, if and only if the form is even (which
means that $b(\xi,\xi)=0$ for any $\xi$).

\pf{Sketch of proof of 1.2.A} Take any class $\xi\in H_n(X;\,\ZZ_2)$,
realize it by a cycle $C$ transversal to $F$. (Here we assume, for
simplicity, that $X$ and $\Gt$ are smooth or piecewise linear.
Lemma is true in general situation, but we need it only for
algebraic varieties.) By a small isotopy of $C$ we can put $C$ into
general position with respect to $\conj(C)$. The intersection
$C\cap\conj(C)$ is a finite set invariant with respect to $\conj$. It
consists of one-point orbits laying in $F$ and two-point orbits
disjoint with $F$. Therefore the number of points of $C\cap\conj(C)$
is congruent modulo 2 to the number of points of $C\cap F$. The first
 of these numbers reduced modulo 2 is $\xi\circ\conj_*\xi$, while
the second one reduced modulo 2 is $\xi\circ\conj_*[F]$ \qed \endpf

Let $A$ be a non-singular
$n$-dimensional real algebraic variety, $\ra$ the set of its real
points and $\ca$ the set of its complex points. Denote by $\conj$ the
complex conjugation involution $\ca\to\ca$. The set $\ra$ is the fixed
point set of $\conj$.

According to Lemma 1.2.A the homology class realized by $\ra$ is
the characteristic class of the $\ZZ_2$-form of $\conj$. It gives the
following theorem.

\tm{1.2.B Theorem} Real algebraic variety $A$ of even dimension
bounds in complexification if and only if the $\ZZ_2$-form of complex
conjugation involution $\conj:\ca\to\ca$ is even.  \qed \endtm

\tm{1.2.C Corollary} A real algebraic surface $A$ of the
3-dimensional projective space having odd degree can not bound in
complexification.
\endtm

\pf{Proof} Take any real plane section $B$ of $A$. The complex point set
$\C B$ of $B$ is invariant under $\conj$, and $\conj$ reverses orientation
of $\C B$. Therefore the class $\Gb\in H_2(\ca;\,\Z_2)$ realized by $\C B$
is invariant under $\conj_*$.
On the other hand it has self-intersection number $\Gb\circ\Gb$ equal to
the degree of $A$ (which is odd by assumption) reduced
modulo 2.  Therefore, the form of involution $\conj$ takes non-zero
value $\Gb\circ\conj_*(\Gb)$ on it. \qed \endpf

These arguments admit straightforward generalization to the situation
of a projective even-dimensional non-singular variety. Remind that the complex
point set $\ca$ of projective variety $A$ of dimension $n$ realizes a non-zero
class belonging to $H_n(\C P^N)=\Z$, this class is equal to d-fold
generator of $H_n(\C P^N)$ (the latter is realized by $\C P^n$). The
number $d$ is called the {\it order  \/} of $A$. It is equal to the
intersection number of $\ca$ and $\C P^{N-n}$ in $\C P^N$. The order of
a hypersurface is its degree, the order of a regular complete
intersection of a collection of hypersurfaces equals the product of the
degrees of these hypersurfaces.

\tm{1.2.D Generalization of 1.2.C} A real projective even-dimensional
non-singular variety of odd order can not bound in complexification.
\endtm

\pf{Proof} Let $A$ be a real non-singular $n$-dimensional subvariety of
order $d$ of the projective $N$-dimensional space. Take a projective real
$(N-\frac n2)$-dimensional subspace $P$ of the ambient space transversal
to $A$.
Denote the intersection of $P$ with $A$ by $B$. Since the self-intersection
of $B$ in $A$ can be obtained as the intersection of $A$ with the
self-intersection of $P$ in the ambient space, the self-intersection number
of $\C B$ in $\ca$ is equal to $d$. Since $\C B$ is invariant under $\conj$,
the class $\Gb\in H_{n/2}(\ca;\,\Z_2)$  is invariant under $\conj_*$. It
has non-zero self-intersection number, since it equals $d\mod2$. Therefore
the form of involution $\conj:\ca\to\ca$ is not even. \qed
\endpf

After these assertions, which mean that under some conditions a real
variety can not bound in complexification, it is natural to wonder, if
it can bound in complexification in any case except the case of curves.
Here is a sufficient condition, providing examples of varieties of any
even dimension which bound in complexification.

\tm{1.2.E Theorem} Any even-dimensional non-singular $M$-variety
with even intersection form of complexification bounds in
complexification.  \endtm

Remind that a real non-singular algebraic variety $A$ is called
{\it $M$-variety,\/} if $\dim H_*(\ra;\,\Z_2)=\dim H_*(\ca;\,\Z_2).$ For
any real algebraic variety $A$ one has
$$\dim H_*(\ra;\,\Z_2)\le\dim H_*(\ca;\,\Z_2).\tag1$$
This inequality was found by R.~Thom \cite{Th}, he observed that it
follows from a general Smith's theorem on homology of fixed
point set of involution (see e.~g. \cite{Br}, Ch. III).
It is a generalization of Harnack's inequality, which says that the number
of components of a real point set $\ra$ of non-singular real algebraic curve
$A$ is at most $g(A)+1$, where $g(A)$ is the genus of $A$.
Inequality (1) is called {\it generalized Harnack inequality.\/}
$M$-varieties are extremal cases of (1). In \cite{V1} I constructed
$M$-surfaces of any degree in 3-dimensional projective space.

\tm{1.2.F Corollary} Any non-singular M-surface of even degree
in projective 3-space bounds in complexification.
\endtm
\pf{Proof} It follows from 1.2.E, since as it is well
known, complexification of any surface $A$ of even degree in 3-dimensional
projective space has even intersection form. \qed
\endpf

1.2.F shows that at least for any even $m>0$ there exists
a real non-singular surface in 3-dimensional projective space of
degree $m$, which bounds in complexification.

\pf{Proof of 1.2.E} It is known (see Rokhlin \cite{R1}) that the
complex conjugation involution of any $M$-variety $A$ acts trivially
in $H_*(\ca;\,\Z_2)$. Therefore the $\Z_2$-form of this involution
coincides with the intersection form of $\ca$. \qed \endpf

\bh 1.3 No more relation among homology classes of real components\endbh
The assertion that a non-singular real algebraic variety $A$ bounds in
classification can be reformulated by saying that the sum of
$\Z_2$-homology classes realized by components of $\ra$ is equal to
zero. It suggests a question, if it can happen that there are other
relations among those classes. In the case of curves the answer is
known to be negative.

Consider the case of $A$ with
$H_1(\ca;\,\Z_2)=0$. The latter equality is not very restrictive
assumption, since all projective regular complete intersections and
cyclic branched covering of projective plane branched over non-singular
plane curve satisfy it.

\tm{1.3.A} If $A$ is a non-singular real algebraic surface $A$ with
$H_1(\ca;\,\Z_2)=0$, then the kernel of the inclusion homomorphism
$H_2(\ra;\,\Z_2)\to H_2(\ca;\,\Z_2)$ has dimension at most 1.
\endtm

\pf{Proof} Consider the Smith sequence of involution
$\conj:\ca\to\ca$.

 $$\minCDarrowwidth{5mm}
\CD
0@>>>H_4(\ca/\conj,\ra)@>>>H_4(\ca)@>>>H_4(\ca/\conj,\ra)@>>>\\
@.          @|             @|                @|            \\
@.      \ZZ_2  @.     \ZZ_2    @.      \ZZ_2    @.          \\
\endCD
$$
$$\minCDarrowwidth{5mm}
\CD
@>>>H_3(\ca/\conj,\ra)@>>>H_3(\ca)@>>>H_3(\ca/\conj,\ra)@>>>\\
    @.  @|                  @|                 @|          \\
    @.   ?      @.           0    @.          \ZZ_2
\endCD $$
$$\minCDarrowwidth{6mm}
\CD
 @>>>H_2(\ca/\conj,\ra)\oplus H_2(\ra)@>>>H_2(\ca) \endCD$$
The group $H_4(\ca)$ is $\Z_2$ since $\ca$ is a closed connected
4-manifold. The  group $H_4(\ca/\conj,\ra)$ is $\Z_2$ by the same
reason (since $\ca/\conj$ is a closed connected manifold). From
that, assumption $H_1(\ca;\,\Z_2)=0$ (which means by
Poincar\'e duality that $H_3(\ca;\,\Z_2)=0$) and exactness of the Smith
sequence it follows that the first boundary homomorphism of it
$H_4(\ca/\conj,\ra)@>>>H_3(\ca/\conj,\ra)$ is bijective, and therefore
$H_3(\ca/\conj,\ra)=\Z_2$. Consequently, the last homomorphism of the
piece of the Smith sequence reproduced above has one-dimensional
kernel. This homomorphism $H_2(\ca/\conj,\ra)\oplus
H_2(\ra)@>>>H_2(\ca)$ restricted to the second summand is the inclusion
homomorphism.  Therefore the kernel of this inclusion homomorphism has
dimension at most 1 (and it happens exactly when the image
of the preceding homomorphism
$H_3(\ca/\conj,\ra)@>>>H_2(\ca/\conj,\ra)\oplus H_2(\ra)$ of the Smith
sequence is contained in the second summand $H_2(\ra)$).\endpf

\ah 2. Complex orientations of a real surface
bounding in complexification\endah

\bh 2.1 Digression: Arnold's principle\endbh   I have to motivate my
further considerations. Otherwise they look more tricky than they
were. Unexpected help came from Arnold's speculations on
complexifications.

Arnold told that everything in mathematics has a complexification.
Some objects have obvious complexification. For example, the
complexification of $\R$ is of course $\C$. Some complexifications are
less obvious. The complexification of the group $\Z_2$ is $\Z$.

Sometimes complexification is not unique. For example, $S^1=\R/\Z$ is a
complexification of $\Z_2$ too. Indeed, $S^1$ is isomorphic to the
multiplicative group of complex numbers with absolute value 1, while
$\Z_2$ is isomorphic to the multiplicative group of real numbers with
absolute value 1. This non-uniqueness of complexification is
explained by the fact that $\Z$ and $S^1$ are dual to each other
(each of them is the character group for the other), while the group
$\Z_2$ is self-dual.

Complexification of symmetric group $S_n$ is the braid group $B_n$.
(Note that this agrees with the previous example: $S_2=\Z_2$ and
$B_2=\Z$!)

Even the most ``real'' notions have complexifications. Consider the
notion of {\it inequality.\/} What is its complexification? To answer to
this question, let us reformulate the statement that $a\ge0$ without
the sign $\ge$. It is easy: assertion $f\ge0$ is equivalent to the
assertion that there exists a real $\xi$ such that $f=\xi^2$.
Inequality is related to the notion of {\it manifold with
boundary\/} (generic inequality defines a manifold with boundary).
In the complex domain equation $f=\xi^2$ is related to the
notion of {\it double branched covering\/}: if the equation
$f(x_1,\dots,x_n)=0$ determines (locally or globally) a complex
hypersurface $A$ of a complex manifold $B$ (where $x_1$,\dots, $x_n$
are coordinates), then the equation $f(x_1,\dots,x_n)=\xi^2$ defines
in an ambient space with coordinates $x_1$,\dots, $x_n$, $\xi$ a
manifold $X$, which two-fold covers $B$ with projection defined by
$(x_1,\dots,x_n,\xi)\mapsto(x_1,\dots,x_n)$.

Thus the complexification of the notion of manifold with boundary is
the notion of two-fold branched covering. Here the boundary is
complexified by the branch locus.

In some situations this transition from a manifold with boundary to
a double branched covering has a concrete form. For example, let $A$
be a real non-singular plane projective curve of an even degree $2k$
defined by an equation $f(x_0,x_1,x_2)=0$. Then the equation
$f(x_0,x_1,x_2)=\xi^2$ defines a non-singular algebraic surface
in a weighted quasi-projective space. The set of its complex points
is a two-fold covering space of $\cpp$ with branching locus $\ca$.
The set of its real points is projected to a part of $\rpp$
defined by inequality $f(x_0,x_1,x_2)\ge0$.

It was this construction that was used by Arnold in his breakthrough
work \cite{A} to get relations between topology of pair $(\rpp,\ra)$
and topology of 4-dimensional manifold $\ca$ with involution
$\conj$.

\bh 2.2 Complexification of the notion of boundary does
work, suggesting definition of complex orientations of surface\endbh
Transition from real curves to real surfaces resembles slightly
complexification. At least all dimensions are multiplied by 2.

In the definition of complex orientations of a real curve the
crucial step was to consider the set of real points of the curve as
the boundary of a half of the set of its complex points.
A happy idea was to look for ``complexification'' of this
construction.

Let $A$ be a non-singular real algebraic surface bounding in
complexification. By definition it means that $\ra$ realizes $0\in
H_2(\ca;\,\Z_2)$. According to Arnold's philosophy, one should find
a double  covering of $\ca$ branched  over $\ra$.

\tm{2.2.A} There exists a double  covering of $\ca$
branched  over $\ra$. This covering is unique up to
equivalence, if $H_1(\ca;\,\Z_2)=0$. \endtm

This assertion follows from the following well-known classification
theorem.

\tm{2.2.B} Let $A$ be a closed $n$-submanifold of an
$(n+2)$-manifold $X$. Double coverings of $X$ branched over $A$
considered up to equivalence are in one to one correspondence with
homology classes $\eta\in H_{n+1}(X,A;\,\Z_2)$ such that
$\p(\eta)=[A]\in H_n(A;\,\Z_2)$. \endtm

I will not prove here 2.2.B, but remind the most classic way to
construct such a covering. Cut the ambient manifold $X$ along a chain
$C$ realizing $\eta$. Take two copies of the result, and glue them
identifying copies of opposite edges of the cut. The result is the
desired covering space.

Deduce 2.2.A from 2.2.B. Consider homology sequence of $(\ca,\ra)$:
$$
H_3(\ca;\,\Z_2)\to
H_3(\ca,\ra;\,\Z_2)@>\p>>H_2(\ra;\,\Z_2)@>{in_*}>>H_2(\ca;\,\Z_2).
$$
By 1.3.A the kernel of $in_*$ contains at most one non-zero element,
and since $A$ bounds in complexification, $\ker(in_*)$ contains
$[\ra]$. By exactness of the homology sequence above, there exists
$\eta\in H_3(\ca,\ra;\,\Z_2)$ with $\p(\eta)=[\ra]\in
H_2(\ra;\,\Z_2)$. If $H_1(\ca;\,\Z_2)=0$, then by Poincar\'e duality
$H_3(\ca;\,\Z_2)=0$, and such $\eta$ is unique.\qed

For the sake of simplicity consider the case $H_1(\ca;\,\Z_2)=0$.
Denote the double covering of $\ca$ branched over $\ra$ by
$DA\to\ca$. Involution $\conj:\ca\to\ca$ preserves $\ra$ and
therefore the induced involution preserves class $\eta\in
H_3(\ca,\ra;\,\Z_2)$ with $\p(\eta)=[\ra]\in H_2(\ra;\,\Z_2)$
characterizing the covering. Consequently $\conj$ can be lifted to
$DA$. In fact there are two liftings of it differing from each other
by the non-trivial automorphism of the covering space $DA$. Denote
them by $c_+$ and $c_-$.

Consider a fiber $D$ of a tubular neighborhood of $\ra$ in
$\ca$. It is homeomorphic to $D^2$ and without loss of generality
one can assume that $\conj (D)=D$ and $\conj$ acts in $D$ as
symmetry with respect to the center $D\cap \ra$ of $D$.  The preimage
$\tilde D$ of $D$ in $DA$ is a fiber of a tubular neighborhood of
$\ra$ in $DA$.  The liftings $c_+$ and $c_-$ of $\conj$ act in
$\tilde D$ as rotations by $\pm \pi/2$, since they cover the
symmetry. Each of them determines an orientation of $\tilde D$,
namely the orientation with respect to which this is rotation by
$\pi/2$ in positive direction. Thus $c_+$ and $c_-$ determine two
orientations of the normal bundle of $\ra$ in $DA$. These
orientations are opposite to each other. Manifold $DA$ is
naturally oriented: there is orientation with respect to which the
projection $DA\to\ca$ has degree $+2$. Therefore the orientations of
normal bundle of $\ra$ in $DA$ above determine orientations of
$\ra$. This pair of opposite orientations is the desired pair of
complex orientations of $A$.

If $\ra$ is not empty, then the liftings $c_+$ and $c_-$ of $\conj$ are
transformations of order 4. Indeed, their squares $c_+^2$, $c_-^2$
cover the transformation $\conj^2=\id$ and are non-trivial, since they
act non-trivially in $\tilde D$.  Therefore they are non-trivial
automorphisms of the covering $DA\to\ca$ and their squares are
identity.

The orbit spaces of them coincide with $\ca/\conj$ and the
projection  $DA\to\ca/\conj$ is the composition of
$DA\to\ca\to\ca/\conj$. It is a cyclic 4-fold covering branched over
$\ra$.

This covering gives another definition of the complex orientations
of $A$. Namely, the 4-fold cyclic covering of $\ca/\conj$ branched
over $A$ gives an element of $H^1(\ca/\conj\sminus \ra;\,\Z_4)=
\hom(H_1(\ca/\conj\sminus \ra),\Z_4)$, which is a characteristic
class of that covering. The dual homology class belongs to
$H_3(\ca/\conj,\ra;\,\Z_4)$.  Its image under the
boundary homomorphism $$H_3(\ca/\conj,\ra;\,\Z_4)\to H_2(\ra;\,\Z_4)$$
is a fundamental class of $\ra$, since the branch index at each
point of $\ra$ is 4. This class is lifted to an
orientation class belonging to $H_2(\ra)$. It is not difficult to
show that it is one of the complex orientations of $A$ defined
above. This construction gives both complex orientations, since the
characteristic class of the 4-fold covering is defined up to sign.

\bh 2.3 Kharlamov's congruence\endbh   Existence of
the complex orientations defined in the preceding section provides
immediately a simple proof of the following theorem, which was
obtained first by Kharlamov \cite{Kh1}.

\tm{2.3.A Kharlamov congruence} If $A$ is a non-singular real
algebraic surface which bounds in its complexification and has
$H_1(\ca;\,\Z_2)=0$, then
$$
\chi(\ra)\equiv 0 \mod8
$$
\endtm

\pf{Proof} The class realized by $\ra$ in $\ca/\conj$ is divisible
by 4, because there exists the cyclic 4-fold covering of $\ca/\conj$
branched over $\ra$. Therefore the self-intersection number of $\ra$
in $\ca/\conj$ is divisible by 16. On the other hand this
self-intersection number is equal to the self-intersection number
of $\ra$ in $\ca$ multiplied by $2$, and the self-intersection
number of $\ra$ in $\ca$ is $-\chi(\ra)$.\qed  \endpf

\bh 2.4 Complex orientations and classes lifted from the
orbit space of complex conjugation\endbh   Let $A$ be a
non-singular real algebraic surface with $H_1(\ca;\,\Z_2)=0$.  Consider
the inverse Hopf homomorphism $$ p^!:H_2(\ca/\conj)\to H_2(\ca) $$
induced by the natural projection $p:\ca\to\ca/\conj$. (Remind that it
can be defined as the composition
$$\CD
H_2(\ca/\conj)@>\text{Poincar\'e duality}>>H^2(\ca/\conj) \\
@. @Vp^*VV\\
H_2(\ca)@<\text{Poincar\'e duality}<<H^2(\ca)
\endCD$$
and geometrically be described as assignment to
the class of a surface transversal to $\ra$ the class of its preimage
under $p$.)

The composition $p_*\circ p^!$ is the multiplication by the degree
of $p$, i.~e. by 2. Since $H_1(\ca;\,\Z_2)=0$, from the universal
coefficient formula it follows that $H_2(\ca)$ has no elements of order 2.
Therefore $p^!:H_2(\ca/\conj)\to H_2(\ca)$ is a monomorphism. Its image
consists of classes invariant under $\conj_*$, but some classes
invariant under $\conj_*$ do not belong to the image. Moreover,
in terms of this image it is possible to give the following
description of the complex orientations.

\tm{2.4.A} Let $A$ be a non-singular real algebraic surface
with $H_1(\ca;\,\Z_2)=0$ which bounds in complexification. The
complex orientations of $A$ are the only orientations such that the
class $\Ga\in H_2(\ca)$ realized by $\ra$ equipped with the
orientation is equal to $2\Gb$ with $\Gb\in p^!(H_2(\ca/\conj))$.
\endtm

\tm{2.4.B Lemma} Multiplication by 2 transforms any class
$\Ga\in H_2(\ca)$ which is invariant under $\conj_*$ into a class
belonging to the image of $p^!$.  \endtm

\pf{Proof of 2.4.B} Note, first, that the situation in homology
with coefficients in $\Z[1/2]$ is simpler. The image of the
inverse Hopf homomorphism
$$
p^!_{\Z[1/2]}:H_2(\ca/\conj;\,\Z[1/2])\to H_2(\ca;\,\Z[1/2])
$$
coincides with the set of classes invariant under $\conj_*$. See
\cite{Br}.

On the other hand, since $H_1(\ca;\,\Z_2)=0$, there is no 2-torsion in
$H_2(\ca)$ and therefore the coefficient homomorphism $H_2(\ca)\to
H_2(\ca;\,\Z[1/2])$ induced by the inclusion $\Z\to\Z[1/2]$ is
a monomorphism.

Since $\ra\ne \empt$, projection $\ca\to\ca/\conj$
induces epimorphism
$$H_1(\ca;\,\Z_2)\to H_1(\ca/conj;\,\Z_2)$$
and therefore $H_1(\ca/\conj;\,\Z_2)=0$. Thus for the same reason as
above, the coefficient homomorphism $H_2(\ca/\conj)\to
H_2(\ca/\conj;\,\Z[1/2])$ is injective. The coefficient
homomorphisms commute with $p_*$ and $p^!$ Therefore in what follows we
may identify elements of $H_2(\ca)$ and $H_2(\ca/\conj)$ with their
images in $H_2(\ca;\,\Z[1/2])$ and $H_2(\ca/\conj;\,\Z[1/2])$ under the
coefficient homomorphisms.

Thus $\Ga$ is an image of some class $\Gb\in
H_2(\ca/\conj;\,\Z[1/2])$, and therefore $2\Ga=p^!(2\Gb)$. However
$2\Gb=p_*\circ p^!(\Gb)=p_*(\Ga)\in p_*H_2(\ca)\subset
H_2(\ca/\conj)$, and $2\Ga\in p^!H_2(\ca/\conj)\subset H_2(\ca)$.
\endpf

\pf{Proof of 2.4.A} By the definition
of complex orientations, the real part $\ra$ equipped with a complex
orientation realizes $0\in H_2(\ca/\conj;\,\Z_4)$. Therefore in
$H_2(\ca/\conj)$ this surface with the complex orientation realizes
class divisible by 4. Denote this class by $\Gg$ and the result of
the division by $\Gd$, so $\Gg =4\Gd$. The class $2\Ga\in H_2(\ca)$
is equal to $p^!(\Gg)$. Therefore $\Ga=2p^!(\Gd)$.

Suppose now that a class $\Ga'\in H_2(\ca)$, which is realized by $\ra$
equipped with some orientation, is equal to $2p^!(\Gd')$. Then
$2\Ga'$ is an image of the class $\Gg'\in H_2(\ca/\conj)$ realized
by $\ra$ with the same orientation. Since $p^!$ is injective,
$\Gg'=4\Gd'$. Divisibility of $\Gg'$ by 4 gives existence of 4-fold
covering of $\ca/\conj$ branched over $\ra$. It is easy to see that
one of generators of its automorpism group defines on $\ra$ the
orientation which gives $\Ga'$. But as it was shown above, such a
covering is unique and the orientation should be one of the complex
ones. \qed
\endpf

The material of this Subsection emerged in a talk with Kharlamov in
the beginning of July 1980.

\bh2.5 Remark on homology description of a real algebraic surface\endbh
In papers on topology of real algebraic K3 surfaces (see
\cite{Kh2}, \cite{Kh3}, \cite{N}, \cite{Kh4}) the topology is
characterized usually by the following homological data: integer
homology $H_2(\ca)$ of the complexification, intersection form
$H_2(\ca)\times H_2(\ca)\to \Z$, and its isometry
$\conj_*:H_2(\ca)\to H_2(\ca)$ induced by $\conj$. Its curious that
these data in the case of K3 surfaces are really sufficient for
description even up to rough projective equivalence (i.~e.
rigid isotopy and projective isomorphism). Thus the image
of $p^!$ does not contain any new information in the case,
and the complex orientations can be restored from the
homology data above. This makes the following questions interesting:

(1) Is that a special property of K3 surfaces, or it can be
generalized to some wider class of real algebraic surfaces?

(2) How are the complex orientations of K3 surfaces restored from
the homology data?

\ah 3. Relative complex orientations of a real surface\endah

\bh 3.1 Types of real algebraic surfaces revised\endbh   In the
case of curve the group $H_1(\ca;\,\Z_2)$ contains only one
naturally distinguished element: zero. Consequently, there are only two
types of real algebraic curves: curves of type I having $[\ra]=0\in
H_1(\ca;\,\Z_2)$ and curves of type II having $[\ra]\ne0\in
H_1(\ca;\,\Z_2)$.

In the case of surface, the group $H_2(\ca;\,\Z_2)$ contains at
least one other naturally distinguished element: the class of
a hyperplane section. As it is well known, it is not zero. Thus
a new type of real algebraic surfaces appears.

A real algebraic surface $A$ is said to be of {\it type $I_{rel}$,\/}
if the $\ra$ is homologous $\mod 2$ to a hyperplane section. A real
algebraic surface, which bounds in complexification, is said to be
of {\it type $I_{abs}$.\/} All other real algebraic surfaces are said
to be of {\it type II.\/}

For some surfaces  one can find other remarkable elements of
$H_2(\ca;\,\Z_2)$. Namely, there may be classes realized by
(complex) algebraic curves distinct from the class of hyperplane
sections.  Phenomena related to the fact that the real part of a
surface can be homologous to complex cycles realized by algebraic
curves deserves special investigation. However a generic surface of
general type has no classes realized by algebraic cycles distinct
from the class of hyperplane section, and therefore I do not feel
necessity to introduce a collection of types finer than one given
above.

\tm{3.1.A Lemma} Real algebraic surface $A$ is of type
$I_{rel}$ if and only if the class of hyperplane section is the
characteristic class of the $\Z_2$-form of the involution
$\conj:\ca\to\ca$. \endtm

\pf{Proof} It follows from Lemma 1.2.A and uniqueness of the
characteristic class. \qed  \endpf

\tm{3.1.B Theorem} Any non-singular M-surface of odd degree in
3-dimensional projective space is of type $I_{rel}$.
\endtm

\pf{Proof} It follows from 3.2.A since, as it is well known, the
class of plane section of a surface of odd degree in 3-dimensional
projective space is the characteristic class of the $\Z_2$-intersection
form, and for an $M$-surface the $\Z_2$-intersection form coincides
with the $\Z_2$-form of complex conjugation involution. \qed\endpf

\bh 3.2 Complex orientations of a real surface modulo a curve\endbh
Now it is natural to expect some analog of complex orientations for
surfaces of type $I_{rel}$. However instead of an analog we find
a generalization.

 Let $A$ be a non-singular real algebraic surface and $C$ a
real algebraic curve on $A$ such that $\C C$ and $\ra$ realize the same
$\Z_2$-homology class.  As above, assume that
$H_3(\ca;\,\Z_2)=0$. (The situation considered in Section 2.2 appears a
special case of this one, if one allows the curve $C$ to be empty.)

Under these assumptions, the construction described below gives two
orientations  of $\ra\sminus\R C$ opposite to each other. They are
called {\it complex orientations of $A$ modulo $C$.\/}

\tm{3.2.A Lemma} There exists a unique two-fold covering of
$\ca$ branched over $\ra\cup\C C$.\endtm

Here the total space of the covering is not a manifold. It has
singularities over the singularities of the branch locus, i.~e. $\R
C=\ra\cap\C C$. By covering branched over $\ra\cup\C C$ we mean a
natural extension of a covering over the complement
$\ca\sminus(\ra\cup\C C)$ of the branched locus such that the covering
can not be extended to a covering of $(\ca\sminus(\ra\cup\C C))\cup pt$
for any point $pt\in\ra\cup\C C $.

\pf{Proof of 3.2.A} We can use a slight generalization of 2.2.B to
the case when the branch locus is a union of two submanifolds.
According to that generalization, coverings under consideration are in
one to one correspondence with homology classes $\eta\in H_3(\ca,
\ra\cup\C C;\,\Z_2)$ with $\p \eta=[\ra]+[\C C]$. Consider the
following segment of the homology sequence of pair $(\ca, \ra\cup\C
C)$:
$$
H_3(\ca;\,\Z_2)\to H_3(\ca, \ra\cup\C C;\,\Z_2)\to
H_2(\ra\cup\C C;\,\Z_2)\to H_2(\ca;\,\Z_2).
$$
Existence of $\eta$
follows from the assumption that $[\ra]+[\C C]$ is mapped by the last
homomorphism to zero. Uniqueness follows from the assumption that
$H_3(\ca;\,\Z_2)=0$.  \qed \endpf

The rest of construction runs as in the absolute case. A reader, who
feels uncomfortable with singular branched covering, may first delete
$\C C$. Singularities would be deleted, and the situation would be
the same as in Section 2.2, but all varieties become non-compact.

The first obvious question on complex orientations of a surface $A$
modulo a curve $C$ is how they are organized in a neighborhood of $\R
C$. One can imagine two opportunities: an orientation of $\ra \sminus\R
C$ can be extendible or not extendible across $\R C$.

\tm{3.2.B} A complex orientation of a surface $A$ modulo a curve
$C$ is not extendible to an orientation of $\ra$.\endtm

\pf{Proof} First, note that locally all objects involved in the
construction of the orientations above are standard. Indeed, consider a
point $pt\in\R C$. Since it is non-singular for both $A$ and
$C$, it has a neighborhood $U$ in $\ca$ such that there exists a
diffeomorphism $h:U\to D^4$ with \roster
\item $h\circ\conj\circ h^{-1}:(x_1,x_2,x_3,x_4)\mapsto(x_1,-x_2,x_3,-x_4)$,
\item $h(U\cap\C C)=\{x\in D^4\,|\,x_3=x_4=0\}$,
\item $h(U\cap\ra)=\{x\in D^4\,|\, x_2=x_4=0\}.$
\endroster Therefore behavior of a complex orientation of $A$ modulo $C$
at $pt$ should be standard too. Now one can trace the construction
above in the model case, but we will use slightly easier indirect
arguments: consider the example with $A$ being projective plane and $C$
a projective line.  By 3.1.B, $\rpp$ and $\C P^1$ are homologous modulo
2 in $\cpp$, and thus projective plane has complex orientations modulo
a projective line.  Since $\rpp$ is not orientable, a complex
orientation of it modulo line can not be extended across $\R C$.
\qed \endpf

\tm{3.2.C Corollary} If $A$ is a non-singular real algebraic
surface with trivial $H_1(\ca;\,\Z_2)$ and $C$  is a non-singular
curve on it such that $[\ra]+[\C C]=0\in H_2(\ca;\,\Z_2)$,
then $\R C$ realizes the element of $H_1(\ra;\,\Z_2)$ dual to the first
Stiefel-Whitney class of $\ra$. \endtm

\pf{Proof} It follows from existence of complex orientation of $A$
modulo $C$ which is an orientation of $\ra \sminus\R C$ and
impossibility to extend it across $\R C$. \qed \endpf

\tm{3.2.D Corollary} The set of real points of any non-singular
real projective surface of even degree and type $I_{rel}$ in
3-dimensional projective space is contractible in $\R P^3$. \endtm

\pf{Proof} Since $\R P^3$ is orientable, and the set of real points
$\ra$ of any surface $A$ of even degree in $\R P^3$ divides $\R P^3$
into two pieces, $\ra$ is orientable as the common boundary of these
two halves of $\R P^3$. Thus any curve in $\ra$ realizing the element
of $H_1(\ra;\,\Z_2)$ dual to the Stiefel-Whitney class of $\ra$ should
bound in $\ra$.

Now consider any non-singular plane section $C$ of $A$.
By 3.2.C  $\R C$ realizes element dual to the Stiefel-Whitney class.
Thus $\R C$ should bound in $\ra$. Therefore the intersection number of
$\R C$ and any loop on $\ra$ is zero. But this number is equal
to the intersection number of the loop and plane in $\R P^3$.
Consequently, any loop on $\ra$ is contractible in $\R P^3$. It follows
that $\ra$ is contractible in $\R P^3$ itself. \qed \endpf

\de{Remark} There exist surfaces of even degree and type $I_{rel}$
contractible in $\R P^3$. The simplest surface of this kind is the
usual sphere. In fact, for it the group $H_2(\ca;\,\Z_2)$ is generated
by the classes of complex conjugate lines in it. A class of plane
section can be realized by union of a line of one family and a line of
the other family. Since lines of the same family are pairwise disjoint
and each line of one family intersects each line of the other family in
one point and transversally, the intersection number of the class of
plane section and the base classes is one. On the other hand,
exactly one line of one family and one line of the other family pass
through each real point of the sphere. Therefore the intersection
number of $[\ra]$ with both base classes is also one. Since any homology
class is characterized by its intersection numbers with the
base classes, $[\ra]$ equals the plane section class and
sphere is of type $I_{rel}$.\endde

\bh 3.3 Semi-orientations\endbh   Pairs of orientations opposite to
each other occur so frequently here, that one feels a necessity to
introduce a term. A pair of orientations opposite to each other will be
called {\it semi-orientation.\/}

The structure contained in a semi-orientation
is equivalent to a construction assigning to a local
orientation at some point $x_1$ a local orientation of
at any other point $x_2$. This assignment should satisfy
the following conditions. \roster
\item Reversion of the local orientation at $x_1$ implies reversion of
the corresponding local orientation at $x_2$.
\item For any three points $x_1$, $x_2$ and $x_3$ the local orientation
at $x_3$ obtained by the construction from a local orientation at
$x_1$ directly coincides with the local orientation obtained from the
same local orientation at $x_1$ in two steps: first, constructing the
corresponding local orientation at $x_2$ and then applying the
construction to the latter local orientation at $x_2$.
\item If $x_1$ and $x_2$ are connected by path, then the construction
coincides with the transfer of local orientation along the path.
\endroster
Any construction satisfying those conditions gives rise to pair of
orientations opposite to each other. Namely, local orientations
obtained by the construction from one local orientation at some point
$x$ constitute an orientation, and starting from the other orientation
at $x$ one get the opposite orientation. It is easy to see, that the
pair of orientations does not depend on the choice of $x$.

On the other hand, any pair of orientations opposite to each other
gives rise to a construction satisfying the conditions above. Namely,
the construction assigning to a local orientation at $x_1$ the local
orientation at $x_2$ such that both local orientations agree with the
same orientation of the given pair of orientations opposite to each
other.

\bh 3.4 Internal definition of complex orientations
(without 4-fold coverings)\endbh   Let $A$ be a non-singular real algebraic
surface with $H_1(\ca;\,\Z_2)=0$ and $C$ be a non-singular real
algebraic curve in $A$ with $[\C C]+[\ra]=0\in H_2(\ca;\,\Z_2)$. We
admit the case of empty $C$.

According to Section 3.2, there is a complex semi-orientation of
$A$ modulo $C$ opposite to each other. In this Section we consider
an alternative construction for it. It is described
in terms of the corresponding assignment of local orientations (see
Section 3.3 above).

Let $x_1$ and $x_2$ be two points of $\ra\sminus\R C$.
Denote by $D_i$ the fiber of tubular neighborhood of $\ra$ over the
point $x_i$ and by $S_i$ the boundary circle of $D_i$. Take a point
$y_i$ in $S_i$.

Since $\dim\ca=4$ and $\dim(\ra\cup\C C)=2$, the space
$\ca\sminus(\ra\cup\C C)$ is connected. Choose a path
$s:I\to\ca\sminus(\ra\cup\C C)$ connecting $y_1$ with $y_2$.

Consider now some local orientations of $\ra$ at $x_1$ and $x_2$. A
local orientation of $\ra$ at $x_i$ defines an orientation of $D_i$
such that the local intersection number of $\ra$ and $D_i$ at $x_i$ is
$+1$.  This orientation of $D_i$ defines an orientation of $S_i$ (since
$S_i=\p D_i$). Let $u_i$ be a path on $S_i$ with $u_i(0)=y_i$ and
$u_i(1)=\conj(y_i)$ which agrees with the orientation of $S_i$.

Consider the loop $s u_2 (\conj\circ s)^{-1}u_1^{-1}$. It is
zero-homologous in $\ca$ modulo 2, since by hypothesis
$H_1(\ca;\,\Z_2)=0$. Therefore the linking number of it with $\ra\cup\C
C$ is well defined (modulo 2).

\tm{3.4.A Lemma} The linking number of the loop $s u_2
(\conj\circ s)^{-1}u_1^{-1}$ and $\ra\cup\C C$ is zero, iff
the local orientations of $\ra$ involved in construction of
$u_i$ agree with the same complex orientation of $A$ modulo
$C$.\endtm

\pf{Proof} Let $Y\to \ca$ be the two-fold covering branched over
$\ra\cup\C C$. Choose a point $\tilde y_1\in Y$ over $y_1$ and
construct a path $v:I\to Y$ with $v(0)=\tilde y_1$ covering $s$. Denote
$v(1)$ by $\tilde y_2$.  Obviously it lies over $y_2$. Construct paths
$\tilde u_i$ covering $u_i$ and starting at $\tilde y_i$. Denote the
point $\tilde u_i(1)$ by $z_i$.

Assume that the local orientations agree with one of two complex
orientation. Then $z_i$ is the image of $\tilde y_i$ under the
transformation $\widetilde\conj:Y\to Y$ covering $\conj:\ca\to\ca$,
and  path $\widetilde\conj\circ v$ joins points $z_1$ and $z_2$.
Therefore path $v\tilde u_2(\widetilde\conj\circ v)^{-1}\tilde u_1^{-1}$
is a loop covering loop  $su_2(\conj\circ s)^{-1}u_1^{-1}$.  It means
that the latter loop has zero linking number with $\ra\cup\C C$.

If the local orientations do not agree with each of the complex
orientations, then reversing local orientation at $x_1$ make them agree
with one of the complex orientations. This changes loop
$su_2(\conj\circ s)^{-1}u_1^{-1}$ by the loop running once along $S_2$,
and make linking number of $su_2(\conj\circ s)^{-1}u_1^{-1}$ with
$\ra\cup\C C$ to be zero. Since linking number of $S_2$ and $\ra\cup\C
C$ is 1, the original linking number of $su_2(\conj\circ
s)^{-1}u_1^{-1}$ with $\ra\cup\C C$ was 1.
\qed \endpf

To define complex semi-orientation of $A$ modulo $C$ in a manner of
Section 3.3, I have to construct for any local orientation of $\ra$ at
$x_1$ a local orientation of $\ra$ at $x_2$ in such a way that the
construction satisfies the conditions of Section 3.3. Lemma 3.4.A
suggests such a construction. For any local orientation of $\ra$ at
$x_1$ one should choose  a local orientation of $\ra$ at $x_2$ such
that the loop provided by the construction above has zero linking
number with $\ra\cup\C C$. It follows from Lemma 3.4.A that this
construction gives the desired complex semi-orientation of $A$ modulo
$C$.

Note that the conditions of Section 3.3 for this construction can be
easily verified independently of Lemma 3.4.A and the original
construction of complex orientations. Thus the construction of this
Section can be used as a base for the whole theory.

\bh 3.5 Orientations modulo changing curve\endbh   In this Section we
study behavior of complex orientation of a surface modulo curve when
the curve moves.

\tm{3.5.A} Let $A$ be a non-singular real algebraic surface with
$H_1(\ca;\,\Z_2)=0$ and $C_1$, $C_2$ be two non-singular real
algebraic curves on $A$ with $\C C_1$ and $\C C_2$ realizing the same
element of $H_2(\ca;\,\Z_2)$ as $\ra$. Then $\R C_1\cup\R C_2$ divides
$\ra$ into two parts (which are the images of sets of real points of
the two-fold coverings of $\ca$ branched over $\C C_1\cup\C C_2$).
Any complex orientation of $A$ modulo $C_1$ and a complex orientation
of $A$ modulo $C_2$ coincide on one of these parts and are opposite on
the other. \endtm

\pf{Proof} Take a point $x_1\in\ra\sminus(\R C_1\cup\R C_2)$.
Reversing, if necessary, the complex orientation of $A$ modulo $C_2$,
one may assume that the complex orientation of $A$ modulo $C_2$
coincides at $x_1$ with the complex orientation of $A$ modulo $C_1$.
Take any point $x_2\in\ra\sminus(\R C_1\cup\R C_2)$. For local
orientations at $x_1$ and $x_2$ induced by the complex orientation of
$A$ modulo $C_1$, apply the construction of Section 3.4 choosing a path
$s$ with $s(I)$ disjoint from $\ra\cup\C C_1\cup\C C_2$. By Lemma
3.4.A, the linking number of loop $su_2(\conj\circ s)^{-1}u_1^{-1}$
with $\ra\cup\C C_1$ is zero. The linking number of the same loop with
$\ra\cup\C C_2$ is zero, iff at $x_2$ the complex orientation of $A$
modulo $C_2$ coincides with the complex orientation of $A$ modulo
$C_1$.  On the other hand, the linking number of the that loop with
$\ra\cup\C C_2$ is equal to the linking number of it with $\C C_1\cup\C
C_2$. The latter depends only on $x_2$.

To complete the proof we have to show that the dependence is as
in 3.5.A, i.~e. that the set of $x_2\in \ra\sminus(\R C_1\cup\R C_2)$
such that the linking number of loop $su_2(\conj\circ s)^{-1}u_1^{-1}$
with $\C C_1\cup\C C_2$ is zero coincides with the image of the set of
real points of the two-fold covering of $\ca$ branched over $\C
C_1\cup\C C_2$.

The covering does exit since $\C C_1$ and $\C C_2$ realize the same
element of $H_2(\ca;\,\Z_2)$. It is unique, since $H_1(\ca;\,\Z_2)=0$.
Denote the covering space by $Z$ and the non-trivial automorphism of
the covering by $\Gt$. Since the complex conjugation involution
$\conj:\ca\to\ca$ preserves $\C C_1\cup\C C_2$, it can be lifted to
$Z$. One can construct a lifting starting at
$x_1\in\ra\sminus(\C C_1\cup\C C_2)$ taking a point  $z_1\in Z$ over
$x_1$ and assuming that $z_1$ is a fixed point for the lifting. After
that the lifting is constructed by continuity in a unique way. There
are two liftings $c_+, c_-:Y\to Y$ obtained from each other by
composition with $\Gt$.  Assume that $c_+(z_1)=z_1$. Since $c_+$ is a
lifting of $\conj$, the fixed point set of $c_+$ is projected into the
fixed point set $\ra$ of $\conj$. To determine, if a point
$x_2\in\ra\sminus(\R C_1\cup\R C_2)$ belongs to the image, one has to
calculate action of $c_+$ in the preimage of $x_2$. Take a point $z_2$
over $x_2$. Connect $z_1$ with $z_2$ by a path $\tilde w$ in the
complement of the preimage of branch locus. This path covers a path $w$
which connects $x_1$ with $x_2$. Then $c_+(z_2)$ is the end point of
the path $c_+\circ\tilde w$.  Therefore $z_2$ is a fixed point of $c_+$
iff paths $\tilde w$ and $c_+\circ\tilde w$ constitute a closed loop.
This is equivalent to the condition that the homotopy class of
loop $w\conj\circ w$ belongs to the group of the covering. This group
consists of homotopy classes of loops in $\ca\sminus(\C C_1\cup\C
C_2)$ unlinked (modulo 2) with $(\C C_1\cup\C C_2)$. Note that for
appropriate choice of $s$ above loops $w\conj\circ w$ and
$su_2(\conj\circ s)^{-1}u_1^{-1}$ are homotopic. \qed\endpf

\bh 3.7 Conversion complex orientations modulo curve into true
orientations\endbh   Results of the preceding Section allow to improve
the constructions of Section 3.2 and 3.4.  While those constructions
give semi-orientation of $\R A\smallsetminus\R C$ which can not be
extended over $\R C$, in this Section we get a semi-orientation of a
two-fold covering space of $\R A$. The role of $C$ will be reduced: the
result depends only on the homology class realized by $\R C$.

First, let me remind some classic purely topological constructions.
With each codimension 1 closed submanifold $Y$ of a manifold $X$ it is
associated a double covering of $X$. This covering can be constructed
in the following way.  One cuts $X$ along $Y$, takes two copies of the
result, and glue them to each other identifying a side of the cut of a
copy with the opposite side of the cut in the other copy. I will denote
the result by $D_YX$.  There is an obvious projection of $D_YX$ onto
$X$, so $D_YX$ is a two-fold covering space of $X$. Various versions of
this construction is used extensively in elementary expositions on
Riemann surfaces.

The result $D_YX$ of the construction above depends only on the
$\Z_2$-homology class realized by $Y$ in $X$. In fact, if $Y'$ is
another submanifold presenting the same $\Z_2$-homology class as $Y$,
then $Y$ and $Y'$ bound together a domain $H\subset X$, and one can
construct a homeomorphism $D_YX\to D_{Y'}X$ which identifies the copies
of $H$ in the copies of $X\sminus Y$ with the copies of $H$ in the
copies of $X\sminus Y'$ with the same numbers, and copies of the
complementary domain  $X\sminus \Cl H$ in the copies of $X\sminus Y$
with the corresponding domains in the copies of $X\sminus Y'$ with
distinct numbers.  The  resulting  homeomorphism  depends  on  the
choice of
$H$. In the case of connected $X$ it can be chosen in two ways and the
corresponding homeomorphisms differs by the non-trivial automorphism of
the covering. In the case of disconnected $X$ a choice should be
done at each component of $X$.

The construction will be applied below to the following situation. Let
$X$ be the set $\R A$ of real points of a non-singular real
algebraic curve $C$ on $A$. If $H_1(\C A; \Z_2)=0$ then the two-fold
covering $D_{\R C}\R A\to\R A$ depends only on homology class $[\C
C]\in H_2(\C A;\Z_2)$ realized by $\C C$ in $\C A$. In fact, for any
curves $C_1$, $C_2$ with $\C C_1$, $\C C_2$ realizing the same element
of $H_2(\C A;\Z_2)$ there exists a distinguished pair of homeomorphisms
$D_{\R C_1}\R A\to D_{\R C_2}\R A$ which differ from each other by the
automorphism of the covering acting non-trivially in each fiber. The
homeomorphisms of this pair are related with two domains $H\subset\R A$
which are the images of sets of real points of the two-fold coverings
of $\C A$ branched over $\R A_1\cup \R A_2$. Cf. Section 3.6. Since
$H_1(\C A;\Z_2)=0$, there is only one two-fold covering of $\C A$
branched over $\R A_1\cup \R A_2$ and the complex conjugation
involution of $\C A$ can be lifted in two ways. Each of these two
liftings defines its own $H$.

Returning to the abstract topological situation above, assume that
there is an orientation
of $X\sminus Y$, which can not be extended across any component of $Y$.
(By the way, it means that $Y$ realizes the $\Z_2$-homology class
Poincar\'e dual to the first Stiefel-Whitney class $w_1(X)$.) Then in
the construction above we get in a natural way a semi-orientation of
$D_YX$. To construct it, one should take the orientation of one of the
copies of the result of cutting induced by the given orientation of
$X\sminus Y$ and take the opposite orientation of the other copy.
Together they induce an orientation of $D_YX$. It is defined up to
sign, since it depends of the choice of the first copy. The
automorphism of the covering $D_YX\to X$ non-trivial over each point of
$X$ reverses these orientations (i.~e. sends them to each other).  This
is also a classical construction known as the construction of the
orientation covering for $X$.  The resulting semi-orientation of $D_YX$
is not changed, if one reverses the original orientation of $X\sminus
Y$, therefore it depends only on the semi-orientation of $X\sminus Y$.

Let $Y'$ be another submanifold presenting the same $\Z_2$-homology
class as $Y$ and $H$ be a domain of $X$ bounded by $Y\cup Y'$. Assume
that $X\smallsetminus Y'$ is oriented in such a way that this
orientation on $H$ coincide with the orientation above of
$X\smallsetminus Y$ and on $X\smallsetminus\Cl H$ is opposite to it.
Then the homeomorphisms $D_{Y'}X\to D_Y X$ defined by $H$ preserve
semi-orientations defined by those orientations of $X\smallsetminus Y'$
and $X\smallsetminus Y$.

Consider now the situation of Section 3.2 and 3.4. Let $A$ be a
non-singular real algebraic surface and $C$ be a non-singular real
algebraic curve on $A$. Let $H_1(\C A;\Z_2)=0$ and $[\C C]+[\R A]=0\in
H_2(\C A;\Z_2)$. Then the constructions of 3.2 and 3.4 give rise to a
complex semi-orientation on $\R A\smallsetminus\R C$ which is not
extendible across $\R C$. The construction above associates with it a
semi-orientation of $D_{\R C}\R A$. If $C'$ is another non-singular
real algebraic curve on $A$ with $\C C'$ realizing the same
$\Z_2$-homology class as $\C C$, then the homeomorphisms $D_{\R C}\R
A\to D_{\R C'}\R A$ constructed as above preserve the semi-orientation.

In some cases the covering space $D_{\R C}\R A$ is the set of real
points of an appropriate algebraic surface. In particular, it happens
in the case of projective real algebraic surface of type $I_{rel}$. If
$A$ is such a surface and $C$ is any curve on $A$ which is a
transversal intersection of $A$ with a real algebraic hypersurface of
odd degree, then $D_{\R C}\R A$ is naturally embedded into the two-fold
covering of the ambient projective space $\R P^N$. This covering space
is sphere $S^N$. The projection is a regular map and therefore  $D_{\R
C}\R A$ is identified with the real part of algebraic variety.

Thus for a projective real algebraic surface of type $I_{rel}$ the
complex semi-orientation modulo hyperplane section corresponds to some
semi-orientation of the real part of another real algebraic surface:
its preimage  under covering $S^N\to\R P^N$.

\ah 4. Conclusion. Survey of some subsequent results\endah

In this section I mention shortly various further developments. I plan
to give a detailed presentation of them elsewhere.

\bh 4.1 Complex orientations of high-dimensional
varieties\endbh  Definition for complex orientations of a surface bounding in
complexification which is given above in Section 3.4 provides not only
opportunity to understand behavior of complex orientations modulo
changing curve, presented in Section 3.5. It suggests a way for
generalizing of the notion of complex orientations in two directions:
for high-dimensional varieties and for high-dimensional analogs of
orientation.

To begin with, consider lower-dimensional case: reformulate the
definition of complex orientations of curves in spirit of Section 3.4.
Let $A$ be a non-singular curve of type $I$. Given points $x_1$,
$x_2\in\R A$, we have to give a criteria for local orientations of $\R
A$ at $x_1$ and $x_2$, if these local orientations agree with the same
complex orientation of $A$. Denote by $y_i$ a point of the boundary of
a tubular neighborhood of $\R A$ in $\C A$ obtained from $x_i$ by a
shift in direction of the normal vector which is $\sqrt{-1}$ times a
tangent vector directed according to the local orientation of $\R A$.
The pair ${y_1, y_2}$ is a 0-cycle $\Z_2$-homologous to zero. If its
linking number with $\R A$ in $\C A$ is zero $(\mod 2)$, then the local
orientations agree with the same complex orientation of $A$, otherwise
the local orientations do not agree with any complex orientation of
$A$.

In high-dimensional case, one may do a similar process. Let $A$ be a
non-singular $n$-dimensional real algebraic variety bounding in
complexification. Take points $x_1$, $x_2\in\R A$ equipped with local
orientations of $\R A$. Assume that these local orientations are
defined by bases $e_1^1,\dots,e_n^1$ and $e_1^2,\dots,e_n^2$ of tangent
spaces of $\R A$ at $x_1$, $x_2$. Fix a tubular neighborhood of $\R A$
in $\C A$ and denote the fiber of it over $x_i$ by $D_i$ and the
boundary of $D_i$ by $S_i$.

Let $y_i$ be the point of $S_i$ obtained by shift of $x_i$ in direction
of $\sqrt{-1}e_1^i$. Choose a path $s$ connecting $y_1$ with $y_2$ in
$\C A\smallsetminus\R A$. It is possible provided $\dim A>1$. In sphere
$S_i$ connect antipodal points $y_i$ and $\conj (y_i)$ by the meridian
passing through the point obtained from $x_i$ by a shift in direction
of $\sqrt{-1}e_2^i$. Those meridians together with $s$ and $\conj\circ
s$ make a 1-cycle $c^1$. If $n=2$, we have to consider the linking
number of this cycle with $\R A$ in $\C A$, as in Section 3.5.
Otherwise this cycle is $\Z_2$-homologous to zero, provided $H_1(\C
A;\Z_2)=0$. Moreover, if $\pi_1(\C A)=0$, it bounds a disk in $\C
A\smallsetminus\R A$. Take this disk, or a chain $s^1$ in $\C
A\smallsetminus\R A$ with $\partial s^2=c^1$. Consider
$s^2+\conj(s^2)$. This is a 2-chain whose boundary consists of great
circles of $S_1$ and $S_2$. Fill this boundary with geodesic disks in
$S_i$ passing through points obtained from $x_i$ by shifts in
directions of $\sqrt{-1}e_3^i$. Denote the resulting 2-cycle by $c^2$.
In the case $n=3$ we have to consider its linking number with $\R A$.
If this linking number is zero, then the local orientations are
announced to agree with the same complex orientation, otherwise the
local orientations do not agree with any complex orientation. If $n>3$,
then the process should be continued.

I do not mean to discuss here conditions under which this construction
gives a well-defined complex orientation. This question is far from
being trivial. It admits, however, a simple solution in the case of
affine varieties. Details will be given elsewhere.

\bh 4.2 $Spin$-structure of a real algebraic surface bounding in
complexification\endbh  Re\-mind that a $Spin$-structure of a manifold
is a reduction of the structure group of its tangent bundle to group
$Spin$.  It can be described in more homological terms as follows:
$Spin$-structure is equivalent to a pair consisting of an orientation
and a $\Z_2$-valued functional defined on the set of framed loops. By a
framed loop in $n$-dimensional manifold $X$  I  mean  here  a  map
$l:S^1\to  X$ equipped with a continuous field of $(n-1)$-frames
assigning to each $t\in S^1$ a sequence of $n-1$ linear independent
tangent vector of $X$ at point $l(t)$. The functional should satisfy
two conditions: first, it should take equal values on homological
 framed loops; second, it takes non-trivial value on any framed loop
with constant map $l:S^1\to X$ and framing defining non-contractible
loop in the space of all $(n-1)$-frames of $T_{l(S^1)}X$. In the case
of two dimensional $X$ the framings are vector fields.

Let $A$ be a non-singular real algebraic surface bounding in
complexification. The following construction gives a $Spin$-structure.
For the orientation of $\R A$ we take one of the complex orientations
defined in Section 2.2 above. Now let $l:S^1\to\R A$ be a loop equipped
with a vector field. Shift $l$ from $\R A$ to $\ca\smallsetminus\R A$
along this vector field multiplied by $\sqrt{-1}$ and take linking
number of the resulting loop with $\R A$. It is easy to check that this
construction gives a functional on the set of framed loops satisfying
the conditions above.

This construction admits generalizations to high-dimensional case
similar to constructions mentioned in Section 4.1.

Let me remind that in the 2-dimensional case $Spin$-structures admit
also a description in terms of $Z_2$-valued quadratic forms on
one-dimensional homology.  Given a closed orientable surface $F$, by a
$\Z_2$-valued quadratic form on $H_1(F;\,\Z_2)$ one means a mapping
$q:H_1(F;\,\Z_2)\to\Z_2$ such that $q(x+y)=q(x)+q(y)+x\circ y$ for any
$x,y\in H_(F;\,\Z_2)$, where $x\circ y$ denotes the intersection number
of classes $x$ and $y$.  There is a one to one correspondence between
$\Z_2$-valued quadratic forms on $H_1(F;\,\Z_2)$  and $Spin$-structures
on $F$ with a fixed orientation of $F$, see e.g. \cite{J}. The
quadratic form corresponding to a $Spin$-structure assings to the
homology class realized  by a collection of disjoint simple closed
loops the number of those loops (modulo 2) plus the sum of values of
the $Spin$-structure on the loops equipped with a vector field
consisting of non-zero vectors tangent to the loops.

Therefore, for any non-singular real algebraic surface $A$ which bounds
in complexification, there is a natural quadratic form
$q:H_1(\ra;\,\Z_2)\to\Z_2$ which assigns to the class represented by a
collection of disjoint simple closed loops $l_1$, \dots, $l_k$ the
number $k$ of those loops (modulo 2) plus the sum of the linking
numbers of $\ra$ and the loops obtained from $l_i$ by a small shift
in the direction of a vector field obtained from a field of vectors
tangent to $l_i$ by mutplication by $\sqrt{-1}$.

Note, that the construction does not involve $\conj$. Therefore it
makes sense in situations when $\conj$ does not exist. For example, in
the situation of a Lagrangian surface homological modulo 2 to zero in
a symplectic 4-manifold. The Spin-structure has been used in study of
Lagrangian tori in $\C^2$, initiated by L.~Polterovich.

\bh 4.3 $Pin^-$-structure of a real algebraic surface of type
$I_{rel}$\endbh  The set of real points of a real algebraic surface of type
$I_{rel}$ may be non-orientable. The obvious modification of the
construction of Section 4.2 gives a $Spin$-structure on the complement
of any hyperplane section of real algebraic surface of type $I_{rel}$.
But this can be essentially improved.

For non-orientable surfaces, there is an analog of $Spin$-structure
which is called
 $Pin^-$-structure. It can be presented as $\Z_4$-valued
quadratic form on one-dimensional $\Z_2$-homology of the surface. Given
a closed surface $F$, by a $\Z_4$-valued quadratic form on
$H_1(F;\,\Z_2)$ one means a mapping $q:H_1(F;\,\Z_2)\to\Z_4$ such that
$q(x+y)=q(x)+q(y)+2(x\circ y)$ for any $x,y \in H_1(F;\,\Z_2)$ where
$x\circ y$ denotes an intersection number (taking value in $\Z_2$ as
above) and 2 denotes the standard inclusion $\Z_2\to\Z_4$.

To define such a form for the set of real points of a non-singular real
algebraic surface $A$ of type $I_{rel}$, consider a collection of
disjoint embedded loops $l_1$, \dots, $l_k$ which presents the homology
class, for which we want to define the value of our quadratic form.
Consider a hyperplane section $C$ with $\R C$ transversal to the loops.
For each of $l_i$ take a non-zero tangent vector field, multiply it by
$\sqrt{-1}$, shift $l_i$ along the result and denote the linking number
of $\ra$ with the loop obtained by $\Gl_i$. Note that $\Gl_i\in \Z_2$
and $2\Gl_i\in\Z_4$.

The value of the quadratic form on the class is equal to
$2\sum^{k}_{i=1}\Gl_i$ plus $2k \mod4$ plus the number of intersection
points of $\R C$ with $\bigcup_{i=1}^{k}l_i$ reduced modulo 4.

One may check that this rule gives a well defined result and that it is
a $\Z_4$-valued quadratic form.

\enddocument